\documentclass[a4paper,10pt]{article}
\usepackage[utf8]{inputenc}
\usepackage{lmodern} 
\usepackage[T1]{fontenc}
\usepackage{textcomp}
\usepackage[scr=boondoxo,frak=boondox,bb=boondox]{mathalfa}
\title{\normalsize{\textbf{ON REAL ALGEBRAIC LINKS IN THE $3$-SPHERE ASSOCIATED WITH MIXED POLYNOMIALS}}}
\author{Raimundo N. Araújo dos Santos  \and  Eder L. Sanchez Quiceno   } 
\usepackage{graphicx}
\usepackage[all]{xy}
\usepackage{mathrsfs}
\usepackage[usenames,dvipsnames]{color}
\usepackage{makeidx}
\usepackage{subcaption}
\usepackage{graphicx}
\usepackage{float}
\usepackage{indentfirst}
\usepackage{amsmath}
\usepackage{amsthm}
\usepackage{graphicx,color}
\usepackage{setspace}
\usepackage{epigraph}
\usepackage{lmodern}
\usepackage{amsfonts,amssymb,amsmath,latexsym}
\usepackage{verbatim}
\usepackage{lscape}
\usepackage{tikz-cd}
\usetikzlibrary{calc,intersections,through,backgrounds}
\usepackage{xfrac}    
\usepackage{pgfplots}
\usepackage{mathtools}
\usepackage{blkarray}
\usepackage{nicematrix}
\usepackage{multirow}
\usepackage{tikz}
\usepackage{faktor}
\usepackage{braids}
\usepackage{titlesec}
\usepackage{lipsum}
\titleformat{\section}
        [block]
        {\large \mdseries\centering}
        {\thesection}
        {1ex}
        {}
        [
        ]%
        \titleformat{\subsection} [block] {\large \mdseries}
        {\thesubsection} {1ex} {}
        [
        ] 
\usepackage[english]{babel}
\usepackage{hyperref}
\hypersetup{
    colorlinks=true,
    linkcolor=blue,
    filecolor=magenta,      
    urlcolor=cyan,
}
\makeindex
\newtheorem{teo}{Theorem}[section]
\newtheorem{corolario}[teo]{Corollary}
\newtheorem{proposition}[teo]{Proposition}

\newtheorem{definition}[teo]{Definition}
\newtheorem{example}[teo]{Example}
\newtheorem{obs}[teo]{Remark}

\newtheorem{conjecture}[teo]{Conjecture}

\newcommand{\C}{\mathbb{C}}       
\newcommand{\R}{\mathbb{R}}       
\newcommand{\N}{\mathbb{N}}       

\newcommand{\re}[1]{\operatorname{Re}(#1)} 
\newcommand{\im}[1]{\operatorname{Im}(#1)}  

\newcommand{\Addresses}{{
  \bigskip
  \footnotesize
    Raimundo~N.~Araújo~dos~Santos, \textsc{Institute of Mathematics and Computer Science (ICMC) University of São Paulo (USP),
   Avenida Trabalhador São-carlense, 400 - Centro
CEP: 13566-590 - São Carlos - SP}\par\nopagebreak
  \textit{E-mail address:} \texttt{rnonato@icmc.usp.br}

  \medskip

  \noindent Eder~L.~Sanchez~Quiceno, \textsc{Institute of Mathematics and Computer Science (ICMC) University of São Paulo (USP),
   Avenida Trabalhador São-carlense, 400 - Centro
CEP: 13566-590 - São Carlos - SP}\par\nopagebreak
  \textit{E-mail address:} \texttt{ederleansanchez@usp.br}

}}
   \makeatletter
\@namedef{subjclassname@2020}{%
  \textup{2020} Mathematics Subject Classification}
\begin{document}
\maketitle
\begin{abstract}
 In this paper we construct new classes of mixed singularities that provide realizations of real algebraic links in the $3$-sphere. Classifications and characterizations of real algebraic links are still open. These new classes of mixed singularities may help to shed light on the Benedetti-Shiota conjecture, which state that any fibered link on the $3$-sphere is a real algebraic link.  
    \vspace{0.2cm} 

   
     \vspace{0.2cm} 
     
   
\end{abstract}



\section{Introduction}
A link $L$ in $S^3$ is fibered if the complement $S^3 \setminus L$ is a surface bundle over the circle $S^1$ whose fiber is the interior of a compact, oriented surface $F$ with $\partial F = L$. It is classically known that fibered links in the sphere $S^3$ were characterized by Stallings in \cite{Stallings1978}. It was shown that $L$ is fibered if and only if the commutator subgroup of the fundamental group of the complement $\pi_1(S^3\setminus L)$ is finitely generated. A very important class of links in $S^3$ are those that arise from some type of singularities of germ of polynomial maps $(\R^4,0)\to (\R^2,0)$, in the following sense: a polynomial map germ $f:(\R^4,0)\to (\R^2,0)$ has a \emph{weakly isolated singularity} at the origin if we have $\Sigma_f\cap V_f=\{0\}$ as a germ of sets, where $\Sigma_f:=\{x \in (\R^4,0): \text{rank}(Jf(x))<2\}$ is the singular set of $f$ and $V_f:=f^{-1}(0)$ both considered as germs of sets. In the particular case where $\Sigma_f=\{0\}$ we say that $f$ has an \emph{isolated singularity} at the origin. According to \cite{Benedetti_Shiota1998}, a link $L$ in $S^3$ is called \emph{weakly algebraic} (resp. \emph{real algebraic}) if there exist a polynomial map germ $f_L : (\R^4, 0) \to (\R^2, 0)$ with a weakly isolated (resp.\ an isolated) singularity at the origin and a positive real number $\rho_0$ such that $L_{f_L,\rho} := V_{f_L} \cap S^3_\rho$ is ambient isotopic to $L$ in $S^3_\rho \approx S^3$ for every $0 < \rho \leq \rho_0$, where $S^3_\rho$ is the sphere of radius $\rho$ centered at the origin in $\R^4$. In this case, such a polynomial map germ $f_L$ is called a \emph{weak realization} (resp. \emph{strong realization}) of $L$. From this perspective  in \cite{Akbulut_King1981} Akbulut and King have proven that every link $L$ in $S^3$ admits a weak realization.  However, they did not provide any sufficient condition for a weak realization to be a strong realization.  
\vspace{0.3cm}

In a more general way, for a germ of complex analytic polynomials $g:(\C^n,0) \to (\C,0)$ with isolated singularity, the classical result due to Milnor \cite{Milnor1968} guarantees that there exists a small radius $\rho_0$ such that for all $0<\rho\leq \rho_0$ the map $$\arg(g):=\dfrac{g}{|g|}:S^{2n-1}_{\rho} \setminus L_{g} \to S^1$$ is a smooth locally trivial fibration. Since the intersection $L_{g}:=V_g\cap S_{\rho}^{2n-1}$ is transverse for all small radii $0<\rho\leq \rho_0$, the manifolds $L_{g}$ for different values of $\rho$ are isotopic. Therefore, the link type of $L_g$ is well-defined and does not depend on the radius $\rho$.  
For a real analytic map germ with an isolated singularity at the origin $f:(\R^n,0)\to (\R^p,0), \ n>p\geq 2$, Milnor also proved that  if $L_{f,\rho}\neq \emptyset$ it is always a fibered link. 

\vspace{0.3cm}

Back to the case where $n=2$ and $g:(\C^2,0)\to (\C,0)$ is a germ of holomorphic functions with isolated singularity\footnote{Consequently analytically equivalent to holomorphic polynomial function.}, the link $L_{g}$ is classically called an \emph{algebraic link} and it was completely characterized as an union of iterated torus knots with some restrictions in their characteristic pairs. For more details see \cite{Brieskorn_Knorrer1986,Zariski1932}.  Since a germ of holomorphic polynomial $f:(\C^2,0) \to (\C,0)$ can be seen as a real polynomial map $(\re{f},\im{f}):(\R^{4},0) \to (\R^2,0)$, then all algebraic links are real algebraic links.

\vspace{0.3cm}

Classification and characterization of real algebraic links are still open and Benedetti and Shiota conjectured that: 
\begin{conjecture}[\cite{Benedetti_Shiota1998}, p. 588]\label{conj}
`` Every fibered link in $S^3$ is the link type of a polynomial map $f$ with isolated singularity at 0, briefly every fibered link is real algebraic."
\end{conjecture}
A partial answer for this conjecture might be assigned to Looijenga in \cite{Looijenga1971}, where as a consequence of a more general construction he proved that given a $NS$-pair $(S^3, L)$, a special connected sum of $NS$-pairs $(S^3, L) \# (S^3, L)$ performed along a non-empty connected component of $L$ produces a new fibered link denoted by $K=L\#L$ which is an \emph{odd fibered link}. It means that, the link $K$ is invariant under the antipodal map $\alpha(x_1, x_2, x_3, x_4) =(-x_1,-x_2,-x_3,-x_4),$ and its fibration map $\phi:S^3\setminus K \to S^1$ satisfies $\phi(\alpha(x))=-\phi(x)$ and hence one can construct a strong realization of the link $K$.
\vspace{0.3cm}

More recently, Bode \cite{Bode2019,Bode2020, Bode2020_arXiv} following a construction made in \cite{Dennis_Bode2019} and \cite{Perron1982}, studied links $L$ in $S^3$ from a braid theoretic point of view and proved that given a \emph{square braid}\footnote{$B=w^2$ for a braid $w$ on $s$ strands.} $B$ on $s$ strands that closes $L$, one can find a mixed polynomial $p:(\C^2,0)\to (\C,0)$ such that the real polynomial map  $p_L:=(\re{p},\im{p}):(\R^{4},0) \to (\R^2,0)$ is a weak realization of the link $L$. Also, in the particular case where the braid $B$ is $P$-fibered (see Definition \ref{defPfibered}), $p_{L}$ can be constructed as a strong realization of $L$. Thus, it is a real algebraic link.
\vspace{0.3cm}

Other strong realizations of fibered links with special properties can be found in the literature. For instance,  Rudolph  \cite{Rudolph1987} constructed mixed polynomials that are strong realizations of the figure-eight knot, Borromean rings and mirror image of any real algebraic link.  Pichon \cite{Pichon2005} constructed strong realizations of the union of iterated torus knots using mixed polynomials of type $f \cdot \bar{g},$ with $f$ and $g$ holomorphic polynomials without common branches. These classes of real algebraic links are examples of fibered links that admits a strong realization as mixed polynomial and provides an interesting topic to be studied.   

\vspace{0.3cm}

In this paper we study singularities of mixed polynomials $f:(\C^2,0)\to (\C,0)$, and we introduce sufficient conditions for $f$ to be a strong realization of the link $L_f$. These conditions define a class of mixed polynomials which contains the classes of strongly non-degenerate mixed polynomials defined by Oka \cite{Oka2010} and the radial semiholomorphic polynomial constructed by Bode \cite{Bode2019}. These conditions are explored for a family of mixed polynomials that are a product of mixed polynomials. The family of such products polynomials is a generalization of the family studied by Pichon \cite{Pichon2005}.
In Section \ref{section2}, we summarize some important properties of mixed polynomials, for instances, non-degeneracy of Newton boundary and radial and polar actions. In Section \ref{section3} we present new criteria, Theorems \ref {semiholomorphiccriterion1} and \ref{teoaxis1}, in the class of mixed polynomials that provide strong realizations of links. Firstly, in Subsection~\ref{3.1} we define, associated with a mixed polynomial $f$, the subsets $\mathcal{A}(f)$, $\mathcal{P}(f)$ and $\mathcal{N}(f)$ of the set of positive weight vectors  $N^+$, and we prove in Proposition~\ref{A0andA} that  $\mathcal{A}(f)=N^+$ is equivalent to the isolatedness of the singularity as a germ of sets at the origin outside of the axes, i.e., $f$ has no critical point in $B_r^4 \cap \C^{*2}$ for some $0<r\ll 1$. Thus we give a condition to guarantee this isolatedness in terms of these subsets.    
\begin{teo}\label{semiholomorphiccriterion1}
Let $f:(\C^2,0) \to (\C,0)$ be a germ of mixed polynomials such that $\mathcal{N}(f) \subset \mathcal{A}(f)$. Then $\mathcal{A}(f)=N^+$. In particular, if $\mathcal{N}(f) \subset \mathcal{P}(f)$ then $\mathcal{A}(f)=N^+$.
\end{teo}
Secondly, in Subsection \ref{3.2} we improve the result above defining two conditions (A-i) and (A-ii) which guarantee the isolatedness at the origin as a germ of sets at the axes.  
 \begin{teo}\label{teoaxis1}
 Let $f:(\C^2,0) \to (\C,0)$ be a germ of mixed polynomials satisfying (A-i) and (A-ii) for an adjusted sequence $P(f)$. Thus, if $\mathcal{A}(f)=N^+$, then $f$ is a strong realization of $L_f$. 
 \end{teo}
  
To finish in Section~\ref{4} we apply our main results and some ingredients given for mixed polynomials to provide several many classes of strong realizations of real algebraic links, for which a description of its topological type can be done without difficulty.  For instance, families of products of mixed polynomials:
 
 \vspace{0.3cm}

 Let $f=p\cdot q: (\C^2,0) \to (\C,0)$ be a product of  two mixed polynomials where $q$ is convenient and polar weighted homogeneous face type, with sequence of positive weight vectors $Q(q):Q_1 \succ Q_2 \succ \cdots \succ Q_{\ell_q}$(see Subsection~\ref{Newtonboundaryofsquarebraids}). Then, set $z=(u, v),\  v=re^{it}$, and consider the face function$$q_{Q_1}(z,\bar{z})=  cr^{n(q)}e^{n_qit},\  c \in \C^*,$$where $(0,n(q))$ is the intersection of the Newton boundary $\Gamma(q)$ with the $v$-axis, where $n_q$ is a nonzero integer with $n(q)-|n_q|$ a non-negative integer. 
 \begin{teo}\label{newrealalgebraiclinks1generalization1}
Given $f=p\cdot q :( \C^2,0)\to (\C,0)$ as above, consider in addition that it is non-degenerate for any $P \in \mathcal{N}(f)$ and that $p$ is semiholomorphic with $deg_u p=s$, and $q$ is polar $u^s$-compatible. If $n_q \in \cap_{P \in \mathcal{N}(f)} \sigma(p_{P})$ and $P_{\ell_p-1} \succ Q_2$ for sequences $P(p)$ and $Q(q)$, then $\mathcal{A}(f)=N^+$. 
\end{teo}
From Theorem~\ref{newrealalgebraiclinks1generalization1} we construct some classes of strong realizations of real algebraic links. In particular, these classes reinforce the idea that the Benedetti-Shiota conjecture \ref{conj} must be true. 

\section{Mixed Polynomial and Newton Boundary}\label{section2}

To know whether or not a germ of mixed polynomial $(\C^2,0) \to (\C,0)$ has an isolated singularity it is not an easy task. Recently combinatorial and analytic techniques have been developed throwing some light on that. Taking as references the works \cite{ Bode2018, Bode2019, Oka2008, Oka2010}, in this section we introduce the basic definitions and notations used in this paper for mixed polynomials. 

\vspace{0.3cm}

We consider a germ of mixed polynomials at the origin $0 \in \C^2$
\begin{equation}\label{mixedpoly}
f(z,\bar{z})=\sum_{j=1}^{l} c_{\nu_j,\mu_j}z^{\nu_j}\bar{z}^{\mu_j}
\end{equation}
where $z = (u,v)$, $ \bar{z}= (\bar{u},\bar{v})$, $z^{\nu_j} = {u}^{\nu_{j,1}}v^{\nu_{j,2}}$ for $\nu_j = (\nu_{j,1} ,\nu_{j,2})$ (respectively $\bar{z}^{\mu_j}=\bar{u}^{\mu_{j,1}}{\bar{v}}^{\mu_{j,2}}$ for $\mu_j=(\mu_{j,1},\mu_{j,2})$).
The zero set  $V_f= f^{-1}(0)$ is called \textit{the mixed hypersurface}. For a fixed variable $x \in \{u,v,\bar{u},\bar{v}\}$ we say that $f$ is a \textit{$x$-semiholomorphic polynomial}  if $f$ does not depend on the variable $\bar{x}$. With this definition we have that  $f$ is a holomorphic polynomial if and only if $f$ is both $u$-semiholomorphic and $v$-semiholomorphic polynomial. Analogously,  $f$ is  anti-holomorphic polynomial if and only if $f$ is both $\bar{u}$-semiholomorphic and $\bar{v}$-semiholomorphic polynomial.
\vspace{0.3cm}

Let $M_{\nu,\mu} = z^{\nu}\bar{z}^{\mu}$ be a mixed monomial where $\nu = (\nu_1,  \nu_2), \  \mu = (\mu_1, \mu_2)$ and let $P =\mathstrut^t(p_1, p_2)$ be a weight vector. We define the \emph{radial degree} of $M_{\nu,\mu}$, $rdeg_P M_{\nu,\mu}$ with respect to $P$  and  the \emph{polar degree} of $M_{\nu,\mu}$, $pdeg_P M_{\nu,\mu}$ with respect to $P$  by $$rdeg_P M_{\nu,\mu} = \sum^2_{j=1} p_j (\nu_j + \mu_j) \text{ and } pdeg_P M_{\nu,\mu} = \sum^2_{j=1} p_j (\nu_j - \mu_j).$$

\begin{definition}
We say $f(z, \bar{z})$ is a radially weighted homogeneous mixed polynomial (for short, radial mixed polynomial) of radial weight type $(P; d_r)$ if 
$$rdeg_P M_{\nu_j,\mu_j} =d_r \text{ for } j=1,2,\dots,l \text{ with } c_{\nu_j,\mu_j}\neq 0$$
 and polar weighted homogeneous of polar weight type $(Q;d_p)$ if 
$$pdeg_Q M_{\nu_j,\mu_j} =d_p \text{ for } j=1,2,\dots,l \text{ with } c_{\nu_j,\mu_j}\neq 0.$$
 If $f(z, \bar{z})$ is  both radially and polar weighted homogeneous of radial weight type $(P;d_r)$ and polar weight type $(Q;d_p)$ then we say that $f(z, \bar{z})$ is a weighted homogeneous mixed  polynomial.
 \vspace{0.2cm}
 
 If $d_p\neq 0$ then $f$ is called a polar non-zero weighted homogeneous polynomial.  Moreover, if for a positive integer $s$ the product $u^s\cdot f(z,\bar{z})$ has non-zero polar degree then $f$ is called of $u^s$-compatible polar non-zero weighted homogeneous polynomial.
\end{definition}

\subsection{Newton boundary and non-degeneracy of mixed polynomials.}\label{subsection3.1}

For a mixed function $f$ as in Eq.~\eqref{mixedpoly} the \textit{Newton polygon} (at the origin) $\Gamma_{+}(f;z,\bar{z})$  is defined by the convex hull of
$$\bigcup_{c_{\nu,\mu}\neq 0}(\nu+\mu)+(\R^{+})^{2}.$$
We call $\Gamma_+(f)$ the Newton polygon of $f$. The \textit{Newton boundary} $\Gamma(f)$ is defined by the union of compact faces of $\Gamma_{+}(f)$. We say that the germ of mixed polynomials $f(z,\bar{z})$ is \emph{convenient} if  $\Gamma(f)$ intersects the $u$-axis and $v$-axis. Let  $(m(f),0)$ and $(0,n(f))$  be the intersection of $\Gamma(f)$ with the  $u$-axis and $v$-axis, respectively.

\vspace{0.3cm}
Given a positive integer vector $P =\mathstrut^t(p_1, p_2)$, we define a linear function $l_P$ on $\Gamma (f)$ via $l_P(\nu) = p_1 \nu_1+p_2 \nu_2 $ for $\nu=(\nu_1,\nu_2) \in \Gamma(f),$ and
let $\Delta (P, f)= \Delta(P)$ be the face where $l_P$ takes its minimal value. We denote the minimal value of $l_P$ by $d(P, f)$. Note that $d(P, f) = \min \{ rdeg_P z^{\nu} \bar{z}^{\mu}| c_{\nu,\mu}\neq 0\}$. 

\vspace{0.3cm}
We define the face function $f_P(z,\bar{z})$ by $$f_{\Delta(P)}(z,\bar{z})=f_P(z,\bar{z})= \sum_{\nu+\mu \in \Delta(P)} c_{\nu,\mu}z^{\nu}\bar{z}^{\mu}.$$
The face function $f_P (z,\bar{z})$ is a radial mixed polynomial of type $(P,d)$ with $d =d(P;f)$.  

\begin{definition}[\cite{Oka2010}]\label{definitionnondegeneracy}
Let $P$ be a positive weight vector. We say that $f$ is non-degenerate for $P$, if the fiber $f_P^{-1}(0)\cap \C^{*2}$ contains no critical point of the mapping $f_P : \C^{*2} \to \C$, where $\C^{*2}=(\C \setminus \{0\})\times (\C \setminus \{0\}) $. We say that $f$ is strongly non-degenerate for $P$ if the mapping $f_P:\C^{*2} \to \C$ has no critical points. 
A mixed polynomial $f$ is called \textit{non-degenerate} (resp.  \textit{strongly non-degenerate}) if $f$ is non-degenerate (resp. strongly non-degenerate) for any positive weight vector $P$. 
\end{definition}
Classes of mixed polynomials where the notions of non-degeneracy and strong non-degeneracy coincide have been studied in  \cite{Blanloe_Oka2015,Inaba_Kawashima_Oka2018,Oka2008,Oka2010,Oka2015}.
\vspace{0.3cm}

 A mixed polynomial $f$ is called \textit{of $u^s$-compatible  polar non-zero weighted homogeneous face type} if for each 1-face $\Delta$, $f_{\Delta}(z,\bar{z})$ is a $u^s$-compatible polar non-zero weighted homogeneous polynomial.  In this case it will be called \textit{polar $u^s$-compatible} for short. 

\vspace{0.3cm}
In this paper we focus on the classes of $u$-semiholomorphic polynomials which by simplicity we call semiholomorphic. We draw attention to the fact that our results can be easily adapted to any $x$-semiholomorphic polynomial for $x \in \{\bar{u},v,\bar{v}\}$.
\vspace{0.3cm}

Consider $p:(\C^2,0)\to (\C,0)$ a semiholomorphic and radially weighted homogeneous polynomial. In this case if we consider the variable $v=re^{it}$, then $p$ can be written as   
 \begin{equation}\label{radrescaling}
p(u,v,\bar{v})= r^{n} g(u/ r^m , e^{it}),
 \end{equation} 
where, $m, n \in \mathbb{Q} $, for an appropriated semiholomorphic function $g:\C \times S^1 \to \C$.
\vspace{0.3cm}

This class of functions as in \eqref{radrescaling} is associated with weak realizations of links that comes from closures of braids. A particular case was proved in \cite{Bode2018,Bode2019} as shortly described below. Before that consider the following definition given in \cite[Definition 1.1]{Bode2020_arXiv}. 
\begin{definition}\label{defPfibered} Let $g:\C \times S^1 \to \C$ be a semiholomorphic function and suppose that $V_g$ is isotopic to the closure of a braid $B$ and let $\arg(g)=\frac{g}{|g|}:\C \times S^1\setminus V_g \to S^1$. We say that the isotopic class of the braid $B$ is $P$-fibered if the projection map $arg(g)$ is a locally trivial fibration.  
 \end{definition}
 
 Let $B$ a braid with $s$-strands and $C$ a fixed component of the closure of $B$. Consider a Fourier parametrization  $\gamma_{C,j}(t)=\left(F_C\left(\frac{t + 2\pi j}{s_C}\right)+i G_C \left(\frac{t + 2\pi j }{s_C}\right), t\right)$, $  j=1,\dots,s_C$, where $s_C$ is the number of strands forming the component $C$. Denote by $\gamma=\{\gamma_{C,j}\}$ the family of parametrized curves indexed by $C's$ and $j's$. Define the function $g_{\gamma}: \mathbb{C} \times S^1 \to \mathbb{C},$  $g_{\gamma}(u,e^{it}) =\prod_{C}\prod_{j=1}^{s_C} (u - F_{C,j}(t)-i G_{C,j}(t)),$ which by definition $V_{g_{\gamma}}$ is the closure of the braid $B$.
 \vspace{0.3cm}
 
 Given an positive integer $k$ consider the semiholomorphic function $p_k(u,v,\bar{v})= r^{2ks} g_{\gamma}(u/ r^{2k} , e^{it})$. It was proved in \cite{Bode2019} that:
 \begin{itemize}
     \item[(i)] There exists positive integer $k$ such that, $p_k$ is a semialgebraic function and has a weakly isolated singularity at the origin. Moreover, $L_{p_k}$ is isotopic to the closure of $B$.
     \item[(ii)] If $B$ is a squared braid then $p_k$ is a radial semiholomorphic polynomial. Thus, $p_k$ is a weak realization the closure of $B$.
     \item[(iii)] If $B$ is both $P$-fibered and squared then $p_k$ is a strong realization of the closure of $B$. 
 \end{itemize}
 \subsection{Newton boundary for  products of mixed polynomials}\label{Newtonboundaryofsquarebraids}\label{generalnewtonboundaryproposition}\label{facefunctionsf}

For a mixed polynomial $f:(\C^2,0)\to (\C,0)$ we know that $\Gamma(f)$ is formed by a finite union of 0-faces and 1-faces. Thus, we can associate to the set of faces of $\Gamma(f)$ a finite sequence of positive weight vectors whose entries are coprimes, which will be called \textit{the sequence of positive weight vector of $f$} (sequence of $f$, for short) as follows:
\vspace{0.3cm}

 \noindent Denote by  $N^+:=\{P=\mathstrut^t(p_1,p_2)| \ p_1,p_2 \in \N \}$ the set of positive vectors. Consider $P= \mathstrut^t(p_1, p_2),Q =\mathstrut^t(q_1, q_2) \in N^+$ then by definition $$P \succ Q \text{ if and only if } \frac{p_1}{p_2} > \frac{q_1}{q_2}.$$ 

For the Newton boundary $\Gamma(f)$ consider the sequence of positive weight vectors $P_1\succ P_2 \succ \cdots \succ P_{\ell_f}$ associated to $f$ satisfying that $\Gamma(f)=\cup_{j=1}^{\ell_f}\Delta(P_j)$ and for all $0$-face $\Delta \in \Gamma(f)$ there exists $j \in \{1,2,\dots,\ell_f\}$ satisfying $\Delta(P_j)=\Delta$. Note that a weight vector $P$ in a sequence of $f$ representing a 0-face can be chosen from the set $[P]:=\{Q\in N^+: \Delta(Q)=\Delta(P)\}$. Note that $f$ is non-degenerate (resp. strongly non-degenerate) if and only if  $f$ is non-degenerate (resp. strongly non-degenerate) for all positive weight vector of a sequence of $f$. A sequence $P_1\succ P_2 \succ \cdots \succ P_{\ell_f}$  is called \textit{minimal} if for any $P’_1\succ P’_2 \succ \cdots \succ P’_{\ell’_f}$ sequence of $f$  then $\ell’_f> \ell_f$. 
\vspace{0.3cm}

Now, consider the product $f=p \cdot q : (\C^2,0)\to (\C,0)$ of two mixed polynomials $p$ and $q$ with sequences $$P(p): P_1\succ P_2 \succ \cdots \succ P_{\ell_p} \text{ and } Q(q):Q_{1} \succ Q_{2} \succ \cdots \succ Q_{\ell_q}.$$ 
 Then the sequence given by $P_{j}$s and $Q_{j}$s provide a sequence of $f$, which may be not minimal. In this case, $f_P=p_P\cdot q_P$ and thus we can see that:  
 
\begin{itemize}

\item[(i)]  If  $P_{\ell_p-1}\succ Q_2 $ then, the face functions of the  $f$ are 
\begin{align*}
f_{P_j}(z,\bar{z})&=p_{P_j}(z,\bar{z})\cdot  q_{Q_1}(z,\bar{z}), \ j=1,2, \dots, \ell_p, \text{ and } \\ f_{Q_j}(z,\bar{z})&=p_{P_{\ell_p}} (z,\bar{z}) \cdot q_{Q_j}(z,\bar{z}), \ j=1,2, \dots, \ell_q.   
\end{align*}
\item[(ii)]  If  $Q_{\ell_q-1}\succ P_2 $ then, the face functions of  $f$ are 
\begin{align*}
f_{P_j}(z,\bar{z})&=p_{P_j}(z,\bar{z})\cdot  q_{Q_{\ell_q}}(z,\bar{z}), \ j=1,2, \dots, \ell_p, \text{ and } \\ f_{Q_j}(z,\bar{z})&=p_{P_{1}} (z,\bar{z}) \cdot q_{Q_j}(z,\bar{z}), \ j=1,2, \dots, \ell_q.   
\end{align*}
\item[(iii)] For a natural number $k$, one can consider $f_k(z,\bar{z})=p_k(u,v,\bar{v})\cdot q(z,\bar{z})$ where $p_{k}$ is a radial semiholomorphic and convenient of type $(2ks+ms,s; d_r)$  defined as
 \begin{equation}\label{radrescalingwithk}
p_k(u,v,\bar{v})= r^{2ks+ms} g(u/ r^{2k+m} , e^{it})
 \end{equation}
 where the sequence of $p_k$ is given by $P_1 \succ P_2=\mathstrut^t(2ks+ms,s) \succ P_3$. Thus, we can prove that $P_2\succ Q_2$ for a large integer $k$, and   the face functions of $f_k$ are  $(f_k)_{P_1}(z,\bar{z})=p_k(0,v,\bar{v}) \cdot  q_{Q_1}(z,\bar{z})$, $ (f_k)_{P_2}(z,\bar{z})=p_k(u,v,\bar{v}) \cdot  q_{Q_1}(z,\bar{z})$, $(f_k)_{P_3}(z,\bar{z})=u^s q_{Q_1}(z,\bar{z})$  and  $(f_k)_{Q_j}(z,\bar{z})=u^s q_{Q_j}(z,\bar{z}).$ 
 \end{itemize}
 
 \section{Mixed Isolated Singularity}\label{section3}
In this section we state our main results concerning to strong realizations of links in the 3-sphere that arise from mixed polynomials. For $f(u,v,\bar{u},\bar{v}):(\C^2,0)\to (\C,0)$ we denote by $f_u$, the partial derivative of $f$ with respect to the variable $u$, analogously we write $f_{\bar{u}}$, $f_v$ and $f_{\bar{v}}$.  
\begin{definition}\label{conjuntosingular}(\cite{Rudolph1987})
Let $f:(\C^2,0)\to (\C,0)$ be a germ of mixed polynomials. Then the singular set of $f$, denoted by $\Sigma_f$ is defined as the solution of 
\begin{equation}\label{criticalpoint}
    \Sigma_f:=
\begin{cases} 
 r_f:=|f_u|^2-|f_{\bar{u}}|^2+|f_v|^2-|f_{\bar{v}}|^2&=0 \\
      s_f:= f_v \bar{f}_{\bar{u}}- \bar{f}_{\bar{v}}f_u=0 
   
   \end{cases}
\end{equation}
\end{definition}
\noindent Note that the Definition \ref{conjuntosingular} is equivalent to the Oka's definition in \cite[Proposition 8]{Oka2010}. In fact, consider $\frac{f_v}{\bar{f}_{\bar{v}}}=\frac{f_u}{\bar{f}_{\bar{u}}}=\lambda$ and substituting $f_v=\lambda \bar{f}_{\bar{v}} $ and $f_u=\lambda \bar{f}_{\bar{u}}$ in $r_f=0$ we get $|\lambda|=1$. Therefore $(\bar{f}_u,\bar{f}_v)=\bar{\lambda}(f_{\bar{u}},f_{\bar{v}})$, with $|\bar{\lambda}|=1$.

\subsection{Mixed singularities outside of the axes}\label{3.1}
Let $f:(\C^2,0)\to (\C,0)$ be a germ of mixed polynomial, $P=\mathstrut^t(p_1,p_2)\in N^+$ and for $(a_1,b_1)\in\C^2$ let $\gamma_{P}(\tau)=(a_1\tau^{p_1}+\text{h.o.t.},b_1\tau^{p_2}+\text{h.o.t.}),\ 0<t<\delta\ll 1$ an analytic curve. We may consider the following set:    
$$\mathcal{A}(f):=\{P \in N^{+} : \text{for any
 $\gamma_{P}(\tau) \in\Sigma_f \Rightarrow a_1=0 \text{ or } b_1=0$} \}.$$

\begin{proposition}\label{A0andA}
Let $f:(\C^2,0)\to (\C,0)$ be a germ of mixed polynomials then $\mathcal{A}(f)=N^+$ if and only if $f$ has no critical points in $B_r^4 \cap \C^{*2}$ for some $0<r\ll 1$.
\end{proposition}
\begin{proof}
($\Rightarrow$) For $P \in N^+ \setminus \mathcal{A}(f)$ by definition there exists $\gamma_P(\tau) \in \Sigma_f$ with $\gamma_P(0)=0$ and $a_1 b_1 \neq 0$. 
\vspace{0.3cm}

($\Leftarrow$) For any sequence of critical points in $B_r^4 \cap \C^{*2}$ by the  Curve Selection Lemma there exists $\gamma_{P}(\tau) \in \Sigma_f$  with $a_1 b_1\neq 0$ for some $P\in N^+$, which contradicts $A(f)=N^+$.  
\end{proof}
\vspace{0.3cm}

We now describe families of mixed polynomials satisfying Proposition \ref{A0andA}.
\vspace{0.3cm}

Let $f(u,v,\bar{v}):(\C^2,0)\to (\C,0)$ be a $u$-semiholomorphic. Then the real Jacobian of $f$ is  
\[
\begin{array}{cccc}
\begin{matrix}
\qquad   \qquad  \quad   \quad (1) \qquad  & \quad  (2) \qquad  & \quad \qquad (3)  \quad \qquad & \qquad  (4) \qquad \qquad 
\end{matrix} \\
J_{\R}f(u,v,\bar{v})= \left( \ \begin{matrix}
 \re{f_u} &  -\im{f_u}  &  \re{f_v+f_{\bar{v}}}&\re{i f_v-i f_{\bar{v}}}  \\
  \im{f_u} &  \re{f_u}  & \im{f_v+f_{\bar{v}}} &  \im{i f_v-i f_{\bar{v}}}   
\end{matrix}\ \right)
\end{array}
\]
For a semiholomorphic polynomial Definition \ref{conjuntosingular} simplifies as $(u,v,\bar{v})$ is a critical point of $f$ if and only if $f_u(u,v,\bar{v})=0$ and $r_f(u,v,\bar{v})=|f_v(u,v,\bar{v})|^2-|f_{\bar{v}}(u,v,\bar{v})|^2=0$. We can see that $r_f(u,v,\bar{v})$ is the determinant of the sub-matrix given by columns $(3)$ and $(4)$ of $J_{\R}f$. 
\vspace{0.3cm}

Again we write $v=re^{it}$ for the second complex variable. Then a critical point of $f$ is given by $(c(r,t),re^{it})$ where $c(r,t)$ is a critical point in variable $u$ of $f(\cdot,re^{it})$, i.e., $f_u(c(r,t),re^{it})=0$, such that $r_f(c(r,t),re^{it})=0$.
\vspace{0.3cm}

 Denote by $\nu(r,t)=f(c(r,t),re^{it})$ the critical values of $f(\cdot,re^{it})$. If $f$ is radial of type $(P;d_r)$ then $$r_f(c(r,t),re^{it})=|f(c(r,t),re^{it})|^2\frac{\partial \arg(\nu(1,t))}{\partial t}$$
 Note that since $f$ is radial the argument of non-zero critical values does not depend on $r>0$. Denote by
\begin{small}
\begin{align*} \sigma(f)&:= \{ m \in \mathbb{Z} : \frac{\partial \arg(\nu(1,t))}{\partial t}+m\neq 0,\ t\in [0,2\pi], \forall \ \nu(1,t) \not\equiv 0 \text{ of } f  \},
\end{align*}
\end{small}
If $f$ is non-degenerate for $P$ then $r_f(c(r,t),re^{it})=0$ if and only if $0 \notin \sigma(f)$.
Analogously, one can consider $f$ above as a $\bar{u}$-semiholomorphic polynomial instead of $u$-semiholomorphic. In such a case we define $\sigma(f,\bar{u}):=-\sigma(\bar{f},u)$. Similarly, one can defines $\sigma(f,v)$ and $\sigma(f,\bar{v})$.

\begin{proposition}\label{radialweightedcondA} 
Let $f:(\C^2,0)\to (\C,0)$ be a germ of radial mixed polynomials of type $(P,d_r)$. Then, $P \in \mathcal{A}(f)$ if and only if $\mathcal{A}(f)=N^+.$ Moreover, if $f$ is semiholomorphic and non-degenerate for $P$ then $0 \in \sigma(f)$ if and only if $\mathcal{A}(f)=N^+.$ 
\end{proposition}
\begin{proof}
The converse of the first statement is  trivial. Suppose that there exists $Q\in N^+ \setminus \mathcal{A}(f)$ and take  $(a,b):=\gamma_Q(\tau_0) \in \Sigma_f\cap \C^{*2}$ for some fix $\tau_0>0$ small enough. Let $P=\mathstrut^t(p_1,p_2)$ then by the radial action $\gamma_P(\tau)=(a\tau^{p_1},b\tau^{p_2})\in \Sigma_f$ and thus $P \notin  \mathcal{A}(f)$, which is a contradiction.
\vspace{0.2cm}

\noindent In particular, if $f$ is semiholomorphic then $f$ is non-degenerate for $P$ and $0 \in \sigma(f)$ if and only if  $B_r^4 \cap \C^{*2} \cap V_f$ and $B_r^4 \cap \C^{*2} \setminus V_f$ have no critical points, respectively, and these two conditions together are equivalents to $\mathcal{A}(f)=N^+.$ 
\end{proof}

 We define the set of weight vectors associated with $f$.
 \begin{align*}
      \mathcal{P}(f)&:=\{P \in N^{+} : f \text{ is strongly non-degenerate for } P  \}.
 \end{align*} 
Let $P(f):P_1 \succ P_2 \succ \cdots \succ P_{\ell_f}$, for short, $P(f):P_j,\ j=1,\dots,\ell_f$ a sequence of $f$. Since non degeneracy is a condition over face functions thus
\begin{align*}
  \mathcal{P}(f)&=\{[P_j],\ j=1,\dots, \ell_f :f \text{ is strongly non-degenerate for } P_j \}.
\end{align*}
Note that for a semiholomorphic polynomial $$ \mathcal{P}(f)=\{[P_j],\  j=1,\dots \ell_f : f \text{ is non-degenerate for } P_j \text{ and } 0 \in \sigma(f_{P_j})\}.$$ Thus we have a characterization of strong non-degeneracy conditions for the case of semiholomorphic polynomials. 
\begin{corolario}
Let $f:(\C^2,0)\to (\C,0)$ be a germ of mixed polynomial and $P(f):P_1 \succ \cdots  \succ P_{\ell_f}$ a sequence of $f$. Then  $f$ is strongly non-degenerate if and only if $P_j \in \mathcal{P}(f),\  j=1,\dots,\ell_f$. Moreover, if $f$ is semiholomorphic and non-degenerate then $0 \in \cap _{j=1}^{\ell_f}\sigma(f_{P_j})$ if and only if $f$ is strongly non-degenerate.   
\end{corolario}
\begin{proof}
Follows directly from comment above and Proposition \ref{radialweightedcondA}. 
\end{proof}
For a mixed polynomial $f$ we say that it is \textit{$\Sigma$-vanished} for a positive weight vector $P$
if $(r_f)_P$ and $(s_f)_P$ have common roots in $\C^{*2}$. Define $$\mathcal{N}(f):=\{P \in N^{+} : f \text{ is $\Sigma$-vanished for } P  \}.$$ 
When $f$ is $u$-semiholomorphic then the $\Sigma$-vanishing of $f$ for a positive weight vector $P$ simplifies to $(f_u)_P$ and $(r_f)_P$ having common roots in $\C^{*2}$, where $r_f=|f_v|^2-|f_{\bar{v}}|^2$.
\begin{example}\label{ex1condA}
For a fixed $\omega\in \C$ consider the semiholomorphic polynomial $$f(u,v,\bar{v}):=p_{\omega}(u,v,\bar{v})(u^2+v^3)$$ where  $p_{\omega}$ is a weak realization of the closure of the braid $\sigma_1^2$ on $2$-strands $p_{\omega}(u,v,\bar{v})=u^2- 2 \omega (v \bar{v} ) u + \omega^2(v\bar{v})^{2}-(v\bar{v})v^2.$
\vspace{0.3cm}

Note that $f$ is a weak realization of $L$ (whose topological type is the same for all $\omega \in \C$) and 
\begin{align*}
f_u(u,v,\bar{v})=&4u^3- 6 \omega (v \bar{v} ) u^2 +( 2\omega^2(v\bar{v})^{2}-2(v\bar{v})v^2+ 2v^3)u- 2 \omega (v \bar{v} )v^3
\end{align*}
\begin{figure}[H]
\begin{subfigure}[a]{.4\textwidth}
    \begin{tikzpicture}[scale=0.8]
\begin{axis}[axis lines=middle,axis equal,yticklabels={0,,2,4,6,8,10}, xticklabels={0,,2,4,6,8,10},domain=-10:10,     xmin=0, xmax=10,
                    ymin=0, ymax=8,
                    samples=1000,
                    axis y line=center,
                    axis x line=center]  
\fill[yellow!90,nearly transparent] (0,50) -- (10,30) -- (30,0) -- (120,0) -- (120,120) -- (0,120) --cycle;
\draw[line width=2pt,red,-stealth,thick,dashed](5,40)--(25,50)node[anchor=south]{$\boldsymbol{P_2}$};
\draw[line width=2pt,red,-stealth,thick,dashed](20,15)--(50,35)node[anchor=south]{$\boldsymbol{P_4}$};
      \filldraw[blue] (0,500) circle (2.3pt)  node[anchor=east,,dashed] {};
    \filldraw[blue] (100,400) circle (2.3pt)  node[anchor=east,,dashed] {};
     \filldraw[blue] (200,200) circle (2.3pt)  node[anchor=east,,dashed] {};
     \filldraw[blue] (100,300) circle (2.3pt)  node[anchor=east,,dashed] {};
     \filldraw[blue] (300,0) circle (2.3pt)  node[anchor=east,,dashed] {};
   \filldraw[black] (80,70) node[anchor=north ] {$\Gamma_+(f_u)$};
   \addplot[mark=*, line width=1.2pt, blue] coordinates{(3,0) (1,3) (0,5)}; 
\end{axis}
\end{tikzpicture}
\caption{$\Gamma(f_u)$.}
\label{boundaryfu}
\end{subfigure}
\hfil \hfil
\begin{subfigure}[a]{.4\textwidth}
\begin{tikzpicture}[scale=0.8]
\begin{axis}[axis lines=middle,axis equal,yticklabels={0,,2,4,6,8,10}, xticklabels={0,,2,4,6,8,10,12},domain=-10:10,     xmin=0, xmax=12,
                    ymin=0, ymax=12,
                    samples=1000,
                    axis y line=center,
                    axis x line=center]  
\fill[yellow!90,nearly transparent] (0,100) -- (10,80) -- (40,30) -- (60,0) -- (140,0) -- (140,120) -- (0,120) --cycle;
\draw[line width=2pt,red,-stealth,thick,dashed](5,90)--(25,100)node[anchor=south]{$\boldsymbol{P_2}$};
\draw[line width=2pt,red,-stealth,thick,dashed](25,55)--(75,85)node[anchor=south]{$\boldsymbol{P_3}$};
\draw[line width=2pt,red,-stealth,thick,dashed](50,15)--(80,35)node[anchor=south]{$\boldsymbol{P_4}$};
      \addplot[mark=*, line width=1.2pt, blue] coordinates{(6,0) (4,3) (1,8) (0,10)}; 
      \filldraw[blue] (0,500) circle (2.3pt)  node[anchor=east,,dashed] {};
    \filldraw[blue] (100,400) circle (2.3pt)  node[anchor=east,,dashed] {};
     \filldraw[blue] (200,200) circle (2.3pt)  node[anchor=east,,dashed] {};
     \filldraw[blue] (100,300) circle (2.3pt)  node[anchor=east,,dashed] {};
   \filldraw[black] (110,105) node[anchor=north ] {$\Gamma_+(r_f)$};
\end{axis}
\end{tikzpicture}
\caption{$\Gamma(r_f)$.}
\label{boundaryfu}\end{subfigure}
\caption{Newton boundary of $f_u$ and $r_{f}$.}
\label{boundaryfu}
\end{figure}
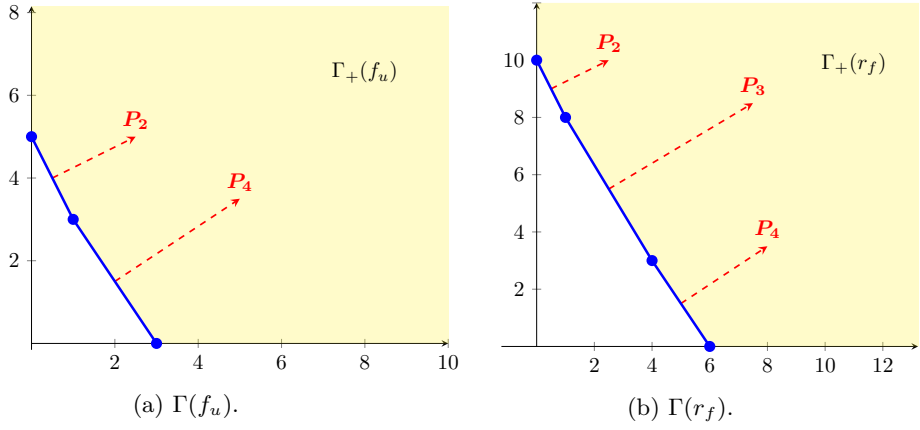
  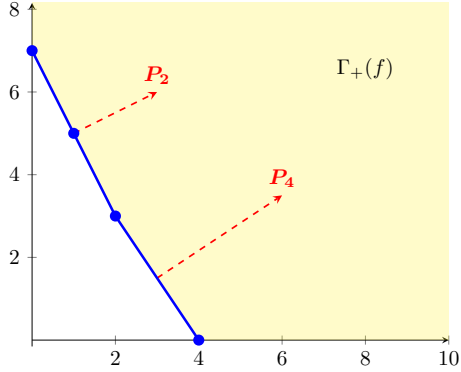
\begin{figure}[H]
\begin{center}
 \begin{tikzpicture}[scale=0.8]
\begin{axis}[axis lines=middle,axis equal,yticklabels={0,,2,4,6,8,10}, xticklabels={0,,2,4,6,8,10},domain=-10:10,     xmin=0, xmax=10,
                    ymin=0, ymax=8,
                    samples=1000,
                    axis y line=center,
                    axis x line=center]   
\fill[yellow!90,nearly transparent] (0,70) -- (20,30) -- (40,0) -- (120,0) -- (120,120) -- (0,120) --cycle;
\draw[line width=2pt,red,-stealth,thick,dashed](10,50)--(30,60)node[anchor=south]{$\boldsymbol{P_2}$};
\draw[line width=2pt,red,-stealth,thick,dashed](30,15)--(60,35)node[anchor=south]{$\boldsymbol{P_4}$};
  \addplot[mark=*, line width=1.2pt, blue] coordinates{(4,0) (2,3) (1,5)(0,7)}; 
  \filldraw[blue] (100,500) circle (2.3pt)  node[anchor=east,,dashed] {};
    \filldraw[blue] (200,400) circle (2.3pt)  node[anchor=east,,dashed] {};
     \filldraw[blue] (300,200) circle (2.3pt)  node[anchor=east,,dashed] {};
  \filldraw[black] (80,70) node[anchor=north ] {$\Gamma_+(f)$};
\end{axis}
\end{tikzpicture}
\end{center}
 \vspace*{-3mm}
\caption{Newton boundary of $f$.}
\label{boundaryfu}
\end{figure}

We want to study the non-$\Sigma$-vanishing of $f$ for a sequence of weight vectors of $f_u$. Let $P(f_u): P_1=\mathstrut^t(3,1)\succ P_2=\mathstrut^t(2,1)\succ P_3=\mathstrut^t(5,3) \succ P_4=\mathstrut^t(3,2)\succ P_5=\mathstrut^t(1,1).$ Note that the face functions of $f_u$ are:
\vspace{0.3cm}

$(f_u)_{P_1}(u,v,\bar{v})=-2 \omega (v \bar{v} )v^3 $,  $(f_u)_{P_3}(u,v,\bar{v})=2v^3u$, $(f_u)_{P_5}(u,v,\bar{v})=4u^3$,
$(f_u)_{P_2}(u,v,\bar{v})=2v^3u -2 \omega (v \bar{v} )v^3$ and  $(f_u)_{P_4}(u,v,\bar{v})=4u^3+2v^3u$.  
\vspace{0.3cm}

Hence, we can see directly that $f$ is non-$\Sigma$-vanished for $P_j, \ j=1,3,5,$. To prove that $f$ is non-$\Sigma$-vanished for $P_2$ and $P_4$, note that $(r_f)_{P_j}= r_{f_{P_j}}$ and $(f_u)_{P_j}=(f_{P_j})_u$ for $j=2,4$. But, the face functions related to $P_2$ and $P_4$ are 
$$f_{P_2}(u,v,\bar{v})=v^3 p_{\omega}(u,v,\bar{v}) \text{ and } f_{P_4}(u,v,\bar{v})=u^2(u^2+v^3).$$ 
After computational calculation we verify that $0 \in \sigma(f_{P_j}), \ j=2,4,$ which implies that $(r_f)_{P_j}$ and $(f_u)_{P_j}$, $j=2,4$ have no common roots in $\C^{*2}$. Therefore, $\mathcal{N}(f)=\emptyset$.
\end{example}
\begin{teo}\label{semiholomorphiccriterion}
Let $f:(\C^2,0) \to (\C,0)$ be a germ of mixed polynomials such that $\mathcal{N}(f) \subset \mathcal{A}(f)$. Then $\mathcal{A}(f)=N^+$. In particular, if $\mathcal{N}(f) \subset \mathcal{P}(f)$ then $\mathcal{A}(f)=N^+$.
\end{teo}
\begin{proof}
 We have to prove $N^+\setminus \mathcal{N}(f) \subseteq \mathcal{A}(f)$. Suppose that there is $P\in N^+\setminus \mathcal{N}(f)$ such that $P \notin \mathcal{A}(f)$ thus we have that a real analytic curve $$\gamma_P(\tau)=(a_1\tau^{p_1}+\text{h.o.t.},b_1\tau^{p_2}+\text{h.o.t.}), \ a_1b_1 \neq 0 , \ 0 \leq \tau \leq 1$$ satisfies that $\gamma_P(\tau)\in \Sigma_f$ for $\tau \neq0$ and $\gamma_P(0)=0$. Since  $\Sigma_f= V_{r_f}\cap V_{s_f}$ thus 
\begin{align*}
   r_f(\gamma_P(\tau))&= (r_f)_P(a_1,b_1)\tau^{d_1}+ \text{h.o.t.}=0, \ d_1=d(P,r_f) \text{ and}\\ 
   s_f(\gamma_P(\tau))&= (s_f)_P(a_1,b_1)\tau^{d_2}+ \text{h.o.t.}=0, \ d_2=d(P,s_f)
\end{align*}
 thus implies that $(r_f)_P(a_1,b_1)=(s_f)_P(a_1,b_1)=0, \ a_1,b_1 \in \C^*$ which is a contradiction to the non-$\Sigma$-vanishing of $f$ for $P$. Therefore $\mathcal{A}(f)=N^+$.
 To finish note that by \cite[Lemma 28]{Oka2010}, $\mathcal{P}(f) \subseteq \mathcal{A}(f)$.  Thus we have the result.
\end{proof}
\begin{example} \label{exemplocloseholomorfo}
\begin{itemize}
    \item[(i)] Consider the semiholomorphic polynomial given by Rudolph's realization of the closure of the braid $\sigma_1 \sigma_2^{-1}$
 $$f(u,v)=u^3-3v \overline{v}u-3v^2\bar{v}u+3v\bar{v}^2u-2 (v+\overline{v}).$$
 We have that $\mathcal{N}(f)=\{(m,m), m \in \N\}$, and by Theorem \ref{semiholomorphiccriterion} it remains to prove that $\mathcal{N}(f)\subset \mathcal{A}(f)$ in order to guarantee  that $f$ has no critical points in $B^{*4}_r$ for some small enough radius $r>0$.
 \item[(ii)] Consider a semiholomorphic polynomial given in example \ref{ex1condA} with $\omega \in \C^*$. By Theorem \ref{semiholomorphiccriterion} we have that there is a small radius $r>0$ such that $B^{*4}_r$ has no mixed singularity of $f$.  
\end{itemize}
\end{example}
\begin{obs}
Note that a mixed polynomial is strongly non-degenerate if and only if $\mathcal{P}(f)=N^+$, which implies triviality that $\mathcal{N}(f)\subset \mathcal{P}(f)$. Therefore, we can see Theorem \ref{semiholomorphiccriterion} as a generalization of Oka's class (\cite[Theorem 19 and Corollary 20]{Oka2010}). For instance, for $\omega=1.1$ of  example~\ref{exemplocloseholomorfo} (ii) belongs to the class of mixed polynomials satisfying $\mathcal{N}(f)\subset \mathcal{P}(f)$ but it is not strongly non-degenerate mixed polynomial. Check for instance $f_{P_1}$. 
\end{obs} 

\subsection{Mixed singularities at the axes}\label{3.2}

 For a mixed polynomial $f$ satisfying $\mathcal{A}(f)=N^+$, we discuss conditions for $f$ to have an isolated singularity at the origin. First, we want to explain a known case of this kind of condition given by Oka in \cite{Oka2010}.  Let $P(f):P_1\succ P_2 \succ \cdots \succ P_{\ell_f} $ a sequence of positive weight vectors such that $P_1, P_{\ell_f}$, are also in the ends of sequences $R(r_f)$ and $S(s_f)$, call this sequences for \textit{adjusted sequence}.
 \begin{proposition}\label{convoka}
 Let $f:(\C^2,0)\to (\C,0)$ be a convenient and strongly non-degenerate mixed polynomial germ. Then $f$ has an isolated singularity.  
 \end{proposition}
 \begin{proof}
 Since $f$ is strongly non-degenerate we have that $\mathcal{P}(f)=N^+$ and by Theorem~\ref{semiholomorphiccriterion} and Proposition~\ref{A0andA} we have that $f$ has no critical points in $\C^{*2}\cap B_r(0)$ for some small radius $r>0$. Therefore, remains to prove that we do not have critical points at the axes close to the origin.  Let $P(f):P_1\succ P_2 \succ \cdots \succ P_{\ell_f} $ be an adjusted sequence. 
   Then, since $f$ is also strongly non-degenerate for $P_1$ and $P_{\ell_f}$ then $f_{P_1}(z,\bar{z})=f_{P_1}(0,v,0,\bar{v})\not\equiv 0 $  and $f_{P_{\ell_f}}(z,\bar{z})=f_{P_{\ell_f}}(u,0,\bar{u},0) \not\equiv  0$ satisfy that the systems $r_{f_{P_j}}=s_{f_{P_j}}=0, \ j=1,\ell_f$ only has trivial solution. But, since $f$ is convenient we have  $s_{f_{P_j}}(z,\bar{z})\equiv 0.$ Hence $r_{f_{P_j}}(z,\bar{z}) \neq 0,\ j=1, \ell_f ,\ z \neq (0,0).$ Consider an analytic curve $\gamma_1(\tau)=(a_1\tau^n +h.o.t,0)\in \Sigma_f$,  
 \begin{equation}\label{eqconvent1}
     (\bar{f}_u(\gamma_1(\tau)),\bar{f}_v(\gamma_1(\tau)))=\lambda(\tau) (f_{\bar{u}}(\gamma_1(\tau)),f_{\bar{v}}(\gamma_1(\tau))), \lambda(\tau)\in S^1
 \end{equation}
  and a Taylor expansion of $\lambda(\tau)=\lambda_0+\lambda_1 \tau+h.o.t.$. We have by convenience of $f$ in $u$ and \eqref{eqconvent1} that $(\overline{f_{P_{\ell_f}}})_u(a_1,0,\bar{a}_1,0)=\lambda_0 (f_{P_{\ell_f}})_{\bar{u}}(a_1,0,\bar{a}_1,0), \ \lambda_0\in S^1$ thus $r_{f_{P_{\ell}}}(a_1,0,\bar{a}_1,0)=0$ which is a contradiction. Briefly,
 $$r_f(\gamma_1(\tau))=s_f(\gamma_1(\tau))=0 \text{ implies }  r_{f_{P_{\ell_f}}}(a_1,0,\bar{a}_1,0)=0,\ a_1\neq 0.$$ Analogously, for $\gamma_2(\tau)=(0, b_1\tau^n +h.o.t)$ we have $$r_f(\gamma_2(\tau))=s_f(\gamma_2(\tau))=0 \text{ implies } r_{f_{P_1}}(0,b_1,0,\bar{b}_1)=0,\ b_1\neq 0.$$
 \end{proof}
 
 Now, we define a generalization of this condition, in order to guarantee the isolatedness of the origin at the axes. Take an adjusted sequence of $f$, $P(f):P_1\succ P_2 \succ \cdots \succ P_{\ell_f} $ and suppose that $f$ satisfies 
 \begin{itemize}
\item[(A-i)] \textit{$(r_f)_{P_{\ell_f}}(u,0,\bar{u},0)=(s_f)_{P_{\ell_f}}(u,0,\bar{u},0)= 0 \text{ if and only if } u=0$} and 
\item[(A-ii)]\textit{$(r_f)_{P_{1}}(0,v,0,\bar{v})=(s_f)_{P_{1}}(0,v,0,\bar{v})= 0 \text{ if and only if } v=0$}. 
\end{itemize}
Thus, following the proof of Proposition~\ref{convoka}, we can prove by contradition.
 \begin{teo}\label{teoaxis}
  Let $f:(\C^2,0) \to (\C,0)$ be a germ of mixed polynomials satisfying (A-i) and (A-ii) for an adjusted sequence $P(f)$. Thus, if $\mathcal{A}(f)=N^+$, then $f$ is a strong realization of $L_f$. 
 \end{teo}
 Note that if $r_f$ and $s_f$ are not convenient in one axis, say $u$, then $f$ does not satisfy (A-i). Moreover, $f$ is not a strong realization of $L_f$. Thus, a necessary condition for $f$ be a strong realization is: for $x=u,v$, $r_f$ or $s_f$ is convenient in $x$. 
 \begin{example}\label{exemplocloseholomorfo2}
In the example \ref{exemplocloseholomorfo} (i) and (ii) we can apply Theorem~\ref{teoaxis} as follow: 
 \begin{itemize}
     \item[(i)]  Consider an adjusted sequence $P(f): P_1=\mathstrut^t(2,1) \succ P_2=\mathstrut^t(1,3)\succ P_3=\mathstrut^t(1,4).$ Note that $f_u$ is convenient in $u$ and $v$ and 
     $(f_u)_{P_{j}}(z,\bar{z})\neq 0, \text{ for all } z \in \C^{*2},\ j=1,3$. Then, by Theorem \ref{teoaxis} $f$ is a strong realization of $L_f$.  
     \item[(ii)] Take $\omega\neq 0$ and an adjusted sequence $P(f):P_1\succ P_2\succ \cdots \succ P_5$. Note that $f_u$ is convenient in $u$ and $v$ and $(f_u)_{P_{j}}(z,\bar{z})\neq 0 \text{ for all } z \in \C^{*2},\ j=1,5$. Then $f_\omega$ is a strong realization of $L_{f_0}$ for any $\omega \in \C^*$. The case $\omega=0$ follows from the convenience and strong non-degeneracy of $f_0$.
 \end{itemize}
\end{example}
\section{Real Algebraic Links}\label{4}
The main application of Theorems \ref{semiholomorphiccriterion} and \ref{teoaxis} is to guarantee that a mixed polynomial is a strong realization of a real algebraic link. It is known that the class of strongly non-degenerate and convenient mixed polynomials has this property. In particular, example \ref{exemplocloseholomorfo2} shows that the class of strong realizations  provided by Theorems \ref{semiholomorphiccriterion} and \ref{teoaxis} is greater than the one given by conditions of strong non-degeneracy and convenience. In this section we show ways to construct explicit examples of strong realizations of real algebraic links. These constructions helps us to reinforce the validity of the Conjecture \ref{conj} and possibly to test invariants of real algebraic links.  
\vspace{0.3cm}

Let $f=p\cdot q: (\C^2,0) \to (\C,0)$ be a product of mixed polynomials where $q$ is convenient and polar weighted homogeneous face type, with sequence $Q(q):Q_1 \succ Q_2 \succ \cdots \succ Q_{\ell_q}$. Then, set $z=(u, v)$, $v=re^{it}$ and consider the face function
$$q_{Q_1}(z,\bar{z})=  cr^{n(q)}e^{n_qit},\  c \in \C^*,$$
where $(0,n(q))$ is the intersection of $\Gamma(q)$ with the $v$-axis, where $n_q$ is a nonzero integer with $n(q)-|n_q|$ a non-negative integer. 
\begin{teo}\label{newrealalgebraiclinks1generalization}
Given $f=p\cdot q :( \C^2,0)\to (\C,0)$ as above, assume in addition that it is non-degenerate for any $P \in \mathcal{N}(f)$ and that $p$ is semiholomorphic with $deg_u p=s$, and $q$ is polar $u^s$-compatible. If $n_q \in \cap_{P \in \mathcal{N}(f)} \sigma(p_{P})$ and $P_{\ell_p-1} \succ Q_2$ for sequences $P(p)$ and $Q(q)$, then  $\mathcal{A}(f)=N^+$. 
\end{teo}
\begin{proof}
 By Theorem \ref{semiholomorphiccriterion} is sufficient to prove that $f$ is strongly non-degenerate for any positive weight vector $P\in \mathcal{N}(f)$. By Section \ref{facefunctionsf} the condition $P_{l_{p}-1} \succ Q_2$ implies that the face functions of $f$ relative to positive weight vectors are: $$f_{Q}(z,\bar{z})=u^s q_{Q}(z,\bar{z})  \text{ and } f_{P}(z,\bar{z})=cr^{n(q)}e^{n_qit} p_{P}(z,\bar{z}),$$ 
where $P$ and $Q$ belong to sequences $P(p)$ and $Q(q)$, respectively. We are going to prove the strong non-degeneracy of $f$ in two parts:
\vspace{0.2cm}

For the first class of face functions $f_Q$, with $Q\in \mathcal{N}(f)$, we have
that $f_Q(z,\bar{z})=u^s q_Q(z,\bar{z})$
is polar non-zero weighted homogeneous polynomial since $q$ is polar $u^s$-compatible. Thus, since $q$ is non-degenerate for $Q$, $f$ is strongly non-degenerate for $Q$.

\vspace{0.2cm}

For the second class of face functions, with $P\in \mathcal{N}(f)$, the non-degeneracy of $f$ for $P$ implies that $f_{P}(z,\bar{z})=cr^{n(q)}e^{n_qit} p_{P}(z,\bar{z})$ has no critical points in $\C^{*2}$ if and only if  for $z=(c(r,t),re^{it}), \ \bar{z}=(0,re^{-it})$ 
\begin{equation}\label{criticalargumentcondition}
    \frac{\partial \arg(f_{P}(z,\bar{z}) )} {\partial t}=n_q+ \frac{\partial \arg(v(r,t))}{ \partial t}\neq 0, \ t\in [0,2\pi], \ r>0,
\end{equation}
where $c(r,t)\not \equiv 0 $ and $v(r,t)=p_P(c(r,t),re^{it},re^{-it})\not \equiv 0$  are respectively,  the nonzero critical roots and nonzero critical values of $p_P(z,\bar{z})$. Now the condition \eqref{criticalargumentcondition} follows from $n_q \in \cap_{P\in \mathcal{N}(f) }\sigma(p_P)$. Thus $f$ is strongly non-degenerate for $P$. Therefore, $\mathcal{A}(f)=N^+$. 
\end{proof}

Note that, in the particular case where $f= p\cdot \bar{q}$ where $p$ and $q$ are holomorphic polynomials without common branches and also convenient and non-degenerate, we can apply our Theorem \ref{newrealalgebraiclinks1generalization} to conclude that $f$ is a strong realization of the link $L_f$. For this reason one can also think that Theorem \ref{newrealalgebraiclinks1generalization} provides a generalization of Pichon's family \cite{Pichon2005} for non-degenerate polynomials.  

\begin{example}
Consider a weak realization of the figure eight knot (see \cite{Bode2018})
\begin{align*}\small
    p(u,v,&\bar{v})=u^3-\frac{3}{4}(v \cdot \overline{v})^{k}u\left(b^2-a^2-ab\left(\frac{v^2}{v \cdot \overline{v}}-\frac{\overline{v}^2}{v \cdot \overline{v}}\right)\right) -\frac{1}{8}(v \cdot \overline{v})^{3k}\cdot \\& \left(a^3\left(\frac{v^2}{v \cdot \overline{v}}+\frac{\overline{v}^2}{v \cdot \overline{v}}\right)+ b^3\left(\frac{v^4}{(v \cdot \overline{v})^2}-\frac{\overline{v}^4}{(v \cdot \overline{v})^2}\right)+3ab^2\left(\frac{v^2}{v \cdot \overline{v}
   }+\frac{\overline{v}
   ^2}{v \cdot \overline{v}}\right)\right)
 \end{align*}
 and the sequences of weight vectors $P(p):P_1\succ P_2 \succ P_3$. Then, for $k$ large and $a\neq b$ we have  $P_{\ell_p-1}\succ Q_2$ and $P \notin \mathcal{N}(p),\ P \succ P_2$. Take a $u^s$-compatible polar germ of convenient semiholomorphic polynomial $q:(\C^2,0)\to(\C,0)$ such that $P_2 \in \mathcal{P}(f)$, i.e., $n_q \in \sigma(p)$. Therefore, by Theorem \ref{newrealalgebraiclinks1generalization} the mixed polynomial $f=p \cdot q :(\C^2,0)\to(\C,0)$ satisfies $\mathcal{A}(f)=N^+$. Now we must apply Theorem \ref{teoaxis} to conclude that $f$ is a strong realization of $L_f$. Note that $f$ is a degenerate mixed polynomial thus $f$ is not in Oka's class of strongly non-degenerate mixed polynomial.     
\end{example}

\begin{corolario}\label{newrealalgebraiclinks1generalization2}
Let $f=p\cdot q :( \C^2,0)\to (\C,0)$ a mixed polynomial and non-degenerate for every $P\in \mathcal{N}(f)$, where $p$ is $x$-semiholomorphic, $x \in \{u,\bar{u}\}$, with $deg_x p=s_p$,  and $q$ $y$-semiholomorphic, $y \in \{v,\bar{v}\}$, with $deg_y p=s_q$. If  $s_q \in \cap_{P \in \mathcal{N}(f)} \sigma(p_{P},x)$ and $s_p \in \cap_{Q \in \mathcal{N}(f)} \sigma(q_{Q},y)$ with $P_{\ell_p-1} \succ Q_2$ for the sequences $P(p)$ and $Q(q)$, then $\mathcal{A}(f)=N^+$.
\end{corolario}
\begin{proof}
For the condition $P_{\ell_p-1} \succ Q_2$ we have that the face functions of $f$ are:  $$f_{Q}(z,\bar{z})=x^{s_p} q_{Q}(z,\bar{z})  \text{ and } f_{P}(z,\bar{z})=y^{s_q} p_{P}(z,\bar{z}),$$ 
where, $P$ and $Q$ belong to sequences of $P(p)$ and $Q(q)$, respectively. Then, considering $s_q \in \cap_{P \in \mathcal{N}(f)} \sigma(p_{P},x)$ and $s_p \in \cap_{Q \in \mathcal{N}(f)} \sigma(q_{Q},y)$, the result follows by  adapting the proof of Theorem \ref{newrealalgebraiclinks1generalization} to this case. 
\end{proof}

Let $p=p_{1,k_1}  p_{2,k_2}$, where  $p_{j,k_j}, \ j=1,2$ are convenient, non-degenerate and radial semiholomorphic polynomials defined as \eqref{radrescalingwithk}, of type $(2k_js_j+m_j s_j,s_j;d_j)$ with $deg_u p_{j,k_j}=s_j$. Note that $p$ is a semiholomorphic polynomial and $deg_u p=s_1 s_2=s$.  Choosing $k_{j's}$ large enough we can consider 
$(2k_1s_1+m_1 s_1,s_1)\succ (2k_2s_2+m_2 s_2,s_2) \succ Q$, for any positive weight vector $Q$. 
\begin{corolario}\label{isolatedsingularitymixedproduct}
Let $f=p\cdot q :( \C^2,0)\to (\C,0)$ be a convenient mixed polynomial and non-degenerate for $P\in \mathcal{N}(f)$, where $p$ is a product of radial semiholomorphic polynomials defined before and $q$ a polar $u^s$-compatible. If $n_q \in \cap_{P\in \mathcal{N}(f)} \sigma(p_{P})$, then for a large enough integers $k_{j's}$ the mixed polynomial  $f=p\cdot q: (\C^2,0)\to (\C,0)$ satisfies $\mathcal{A}(f)=N^+$.
\end{corolario}
\begin{proof}
Since $p$ is a product of radial semiholomorphic polynomials $p_{j,k_j},\ j=1,2$ of type $(2k_js_j+m_js_j,s_j;d_j)$ we have a sequence of positive weight vectors $P(p)$ and $P_{\ell_p-1}=(2k_2s_2+m_2 s_2,s_2)$.  Note that we can choose $k_{j}$s such that the sequence of positive weight vectors $Q(q)$ satisfies $P_{\ell_p-1}\succ Q_2$.  Then the result follows by Theorem~\ref{newrealalgebraiclinks1generalization}. 
\end{proof}
\begin{obs}\label{remark35}
\begin{itemize}
\item[(i)] The above result can be generalized to a family of type $p=p_{1,k_1}  p_{2,k_2} \ldots p_{l,k_l}  $ by making the appropriated modifications. 
\item[(ii)] For a convenient radial semiholomorphic polynomial $p_{k,\omega}$ of type $(2ks+ms,s;d_r)$, with $deg_u p_{k,\omega}=s$ given by
$$p_{k,\omega}(u,v,\bar{v}):=r^{2ks+ms}g\left(\frac{u}{r^{2k+m}}+\omega,e^{it}\right)$$
that is also a weak realization of the closure of a square braid, we can prove that for every complex number $\omega$ with large module we have that for all  $v \in \C^*$  $$-1<\frac{\partial \arg(p_{k,\omega}(0, v, 0, \bar{v}))}{\partial t}<1.$$ Hence, $\cap_{j=1}^3\sigma((p_{k,\omega})_{P_j})=\sigma(p_{k,\omega})\setminus \{0\}$. Moreover, $p_{k,\omega}(z,\bar{z})$ is non-degenerate. 
\end{itemize}
\end{obs}
With these notations and definitions we can prove the following result: 
\begin{proposition}\label{newrealalgebraiclinks1generalization3}
Let  $p_{k, \omega}$ as in the Remark \ref{remark35}, item (ii), and $q:(\C^2,0)\to (\C,0)$ a polar $u^s$-compatible non-degenerate and convenient. If $n_q \in \sigma(p_{k,\omega})$ then there exist an integer $k$ and a complex number $\omega$ such that $f_{k,\omega}=p_{k,\omega}\cdot q: (\C^2,0)\to (\C,0)$ is a strong realization of $L_{f_{k,\omega}}$. Moreover, if $|\omega|$ is large enough then $f$ is also strongly non-degenerate.
\end{proposition}
\begin{proof}
Taking $k>0$ large enough we can consider $P(f):P_1\succ P_2\succ Q_2 \succ \cdots \succ Q_{\ell_q}$. Since $q$ is polar $u^s$-compatible and $n_q \in \sigma(p_{k,\omega})$ then for all $P \preceq P_2$, $P \in \mathcal{P}(f_{k,\omega})$. Thus, by Corollary \ref{isolatedsingularitymixedproduct} we have $\mathcal{A}(f)=N^+$ since for all $P \succ P_2$, we can prove that there exists $\omega$ such that \begin{align}\label{cond5prop3.17}
  (r_{f_{k,\omega}})_P(0,v,0,\bar{v})\neq 0 \text{ or }(s_{f_{k,\omega}})_P(0,v,0,\bar{v}) \neq 0 \text{ for all } v\in \C^{*},   
\end{align}
therefore  $P \notin \mathcal{N}(f_{k,\omega})$. The isolatedness of the singularity at the origin follows from \eqref{cond5prop3.17} and that for all $Q \preceq Q_{\ell_q-1}$,   $$(r_{f_{k,\omega}})_Q(u,0,\bar{u},0)\neq 0 \text{ or } (s_{f_{k,\omega}})_Q(u,0,\bar{u},0)\neq 0 \text{ for all } u\in \C^{*}.$$ Moreover if $\omega$ is large enough by the Remark \ref{remark35}, item (ii) we can  guarantee that for all $P \preceq P_2,$ $P \in \mathcal{P}(f_{k,\omega})$. Thus  $\mathcal{P}(f_{k,\omega})=N^+$ and the strongly non-degeneracy conditions follows. 
\end{proof}
\begin{obs}
Note that $p_{k,\omega}$ is a mixed polynomial for all $\omega \in \C$ only in the case where it is associated to a square braid. However, Proposition \ref{newrealalgebraiclinks1generalization3} can be formulated with a little more extra details for a radial semiholomorphic polynomial $p_k$ instead of $p_{k,\omega}$.   
\end{obs}
\begin{example}\label{exemploprodutostrongly}
 Consider $q : (\mathbb{C}^2,0) \to (\mathbb{C},0)$ given by $q(z,\bar{z})=u^{2}+v^{n}$ and $p_{k}$ a weak realization of the closure of $B=(\sigma_1^{-1} \sigma_2\sigma_1^3\sigma_2)^2$ (closure of square of $B$, figure \ref{figsquare52:sub-second}), constructed as in \cite{Bode2019,Dennis_Bode2017}. Consider a Fourier parametrization of the braid $5_2$ given for the trigonometric functions $$(F'_{5_2}(t),G'_{5_2}(t))=\left(\cos(2t) + \frac{3}{4}\sin(5t), -\sin(4t) - \frac{1}{2}\cos(t)\right).$$  
 Any phase shift of a Fourier parametrization is also a Fourier parametrization. Hence, we consider the phase shift 
 $(F_{5_2}(t),G_{5_2}(t))=(F'_{5_2}(t+\pi/6),G'_{5_2}(t+\pi/6))$. 
 \begin{figure}[h]
  \hfil
\begin{subfigure}[b]{.4\textwidth}
  \begin{center}
  \includegraphics[height=0.5cm,angle=0,width=0.8\linewidth]{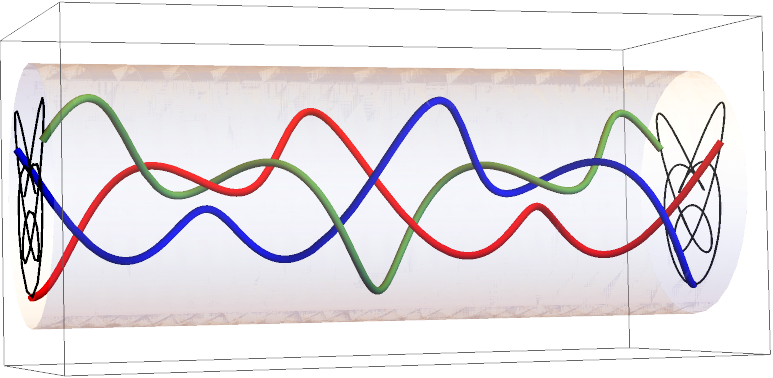}
  \vspace{0.4cm}
  \caption{Braid $(\sigma_1^{-1} \sigma_2\sigma_1^3\sigma_2)^2$.}
  \label{figsquare52:sub-first}
  \end{center}
\end{subfigure}
\begin{subfigure}[b]{.4\textwidth}
  \begin{center}
  \includegraphics[height=0.5cm,angle=0,width=0.8\linewidth]{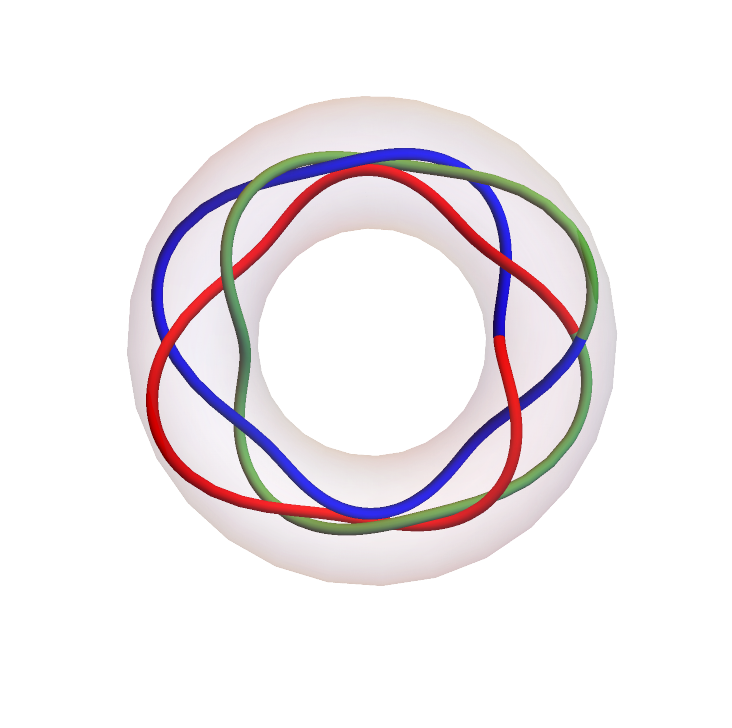}
   \vspace{-0.7cm}
  \caption{Closure of $B$.}
  \label{figsquare52:sub-second}
   \end{center}
\end{subfigure}
  \caption{Square braid formed by the Fourier parametrization of $(F_B(\frac{t+2\pi j}{3}),G_B(\frac{t+2\pi j}{3}),t) \ j=1,2,3 $.}
\label{fig:fig}
\end{figure}
 Let $B$ be the square of a braid, whose closure is the knot $5_2$ (see figure \ref{figsquare52:sub-first}) and  the natural Fourier parametrization of $B$ given by  $(F_{B}(t),G_{B}(t))=(F_{5_2}(2t),G_{5_2}(2t)).$ Thus the mapping $g$ is given by 
$$  g(u,e^{it})=\prod_{j=1}^3 \left( u- \left(F_B \left(\frac{t + 2\pi j}{3}\right)+i \left (G_B\left(\frac{t + 2\pi j }{3}\right) \right) \right)\right).$$
 Expanding the product we have the function 
 \begin{align*}
 g(u,e^{it})=u^3+u \bigg(-\frac{15}{64} + \frac{3}{16} e^{-2 i t}  - \frac{9}{16} e^{2 i t} +\frac{33}{32} e^{-4 i t} - \frac{33}{32} e^{4 i t}  & \\  +\frac{9}{16}  e^{-6 i t} -\frac{9}{16}   e^{6it} \bigg) +\bigg(\frac{187}{128} e^{-2it} - \frac{1}{128}  e^{2it}  &\\+ \frac{115}{128} e^{-4 i t}+ \frac{13}{128} e^{4 i t}+ \frac{51}{256} e^{-6 i t}+\frac{93}{256} e^{6 i t}+\frac{3}{4}   &\\ +\frac{11}{128} e^{-8 i t}+ \frac{43}{128} e^{8 i t} + 
 \frac{27}{512}  e^{-10 i t} + \frac{27}{512}  e^{10 i t}&\bigg).
  \end{align*}

Considering the mixed polynomial $p_{k,\omega}$ given by 
  \begin{align*}
 p_{k,\omega}(z,\bar{z})=(v\bar{v})^{3k}g(u/(v\bar{v})^{k}+w,v/\sqrt{v\bar{v}}),
 \end{align*} 
  then for any $\omega$ we get that $\sigma(p_{k,\omega})=\mathbb{Z} \setminus \{\mathbb{Z} \cap (-24,4)\} $. Let $f_{k,\omega,n}(z,\bar{z})=p_{k,\omega}(z,\bar{z})\cdot (u^2+v^{n})$, we can see that exists $\omega$ such that $(f_{k,\omega,n})_u(0,v,\bar{v})\neq 0$, for all $v\in \C^*$. Hence, since $n>3$ satisfies $n \in \sigma(p_{k,\omega})$ then by Proposition \ref{newrealalgebraiclinks1generalization3} we can choice a large integer $k$ such that the mixed polynomial $f_{k,\omega,n}(z,\bar{z})$ is a strong realization of $L_{f_{k,\omega,n}}$. 
\vspace{0.3cm}

For instance, for $n=5$ the figure \ref{Lflinkex} represents the braid on $5$-strands that closes to a link isotopic to $L_{f_{k, 5,5}}$. 
  \begin{figure}[h]
  \centering
\includegraphics[width=0.250\linewidth, angle=0,height=1.4in]{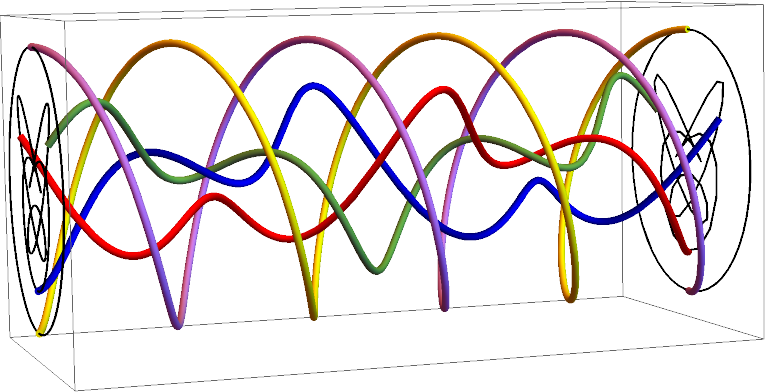} 
\caption{Braid representation of  $L_{f_{k,5,5}}$ with $(u,t) \in \C \times [0,2\pi]$ when $v=r_{0} e^{it}, \text{ and } r_0 \ll 1$.}
 \label{Lflinkex}
\end{figure}

An analogous argument can be done for $h_{k,\omega,n}(z,\bar{z})=p_{k,\omega}(z,\bar{z})\cdot (\bar{u}^5+\bar{v}^n)$. Then since $n>24$ satisfies $-n \in \sigma(p_{k,\omega})$ then by Proposition \ref{newrealalgebraiclinks1generalization3} we can choice a large integer $k$ such that the mixed polynomial $h_{k,\omega,n}(z,\bar{z})$ is strong realization of $L_{k,\omega,n}$ when $\omega$ satisfies $r_{h_{k,\omega,n}}(0,v,\bar{v})\neq 0$ or $s_{h_{k,\omega,n}}(0,v,\bar{v})\neq 0$, for all $v\in \C^*$. To illustrate the case consider $h_{5,0,25}$, then we have that  $\mathcal{N}(h_{5,0,25})=\{[P_2],[P_4]\}$ where 
$P_2=\mathstrut^t(10,1)$ and $P_4=\mathstrut^t(5,1)$. Calculating the face functions of positive weight vectors belonging to $\mathcal{N}(h_{5,0,25})$ we have $(h_{5,0,25})_{P_2}(z,\bar{z})=u^3 (\bar{u}^5+\bar{v}^n)$ and $(h_{5,0,25})_{P_4}(z,\bar{z})=\bar{v}^{25} p_{5,0}(z,\bar{z})$. It is clear that $\mathcal{N}(h_{5,0,25}) \subset \mathcal{P}(h_{5,0,25})$ since $q$ is polar $u^3$-compatible and $-25 \in \sigma(p_{k,0})$.
\end{example} 

Let $p, q  :( \C^2,0)\to (\C,0)$ be germs of non-degenerate mixed polynomials with $p$ convenient and semiholomorphic polynomial and $q$ an germ of holomorphic isolated singularity. Suppose that $P_{\ell_p-1} \succ Q_2$ for sequences $P(p)$ and $Q(q)$ of $p$ and $q$, respectively. 
\begin{proposition}\label{corollarystrongandintrictopologicalinv}
If there exists $c\neq 0$ such that either
\begin{itemize}
\item[(i)]  $q_{Q_1}(u,v)=c u \cdot v^{n}$, 
 and $n \in \cap_{P \in \mathcal{N}(f)} \sigma(u\cdot  p_{P_j})$, or 
 \item [(ii)] $q_{Q_1}(u,v)= c v^{n}$ and $n \in \cap_{P \in \mathcal{N}(f)} \sigma(p_{P_j})$, 
\end{itemize}
then $f=p \cdot q:(\C^2,0)\to (\C,0)$ is a strong realization of $L_f$.
\end{proposition}
\begin{proof}
The condition $P_{l_{p}-1} \succ Q_2$ implies that the face functions of $f$ are: $$f_{Q}(z,\bar{z})=u^{m(p)} q_{Q}(z,\bar{z})  \text{ and } f_{P}(z,\bar{z})=q_{Q_1}(u,v)\cdot p_{P}(z,\bar{z}),$$ 
where $P$ and $Q$ are positive weight vector of the sequences $P(p)$ and $Q(q)$, respectively. Note that a brief calculation shows that $f$ is non-degenerate semiholomorphic polynomial, moreover $f$ is a weak realization of $L_f$. By (i) and (ii) and that the face functions $f_Q$ are holomorphic non-degenerate polynomials we have $\mathcal{N}(f) \subset \mathcal{P}(f)$ thus for Theorem \ref{semiholomorphiccriterion} follows $\mathcal{A}(f)=N^+$. In both cases the non-degeneracy and convenience of $p$ and the non-degeneracy of $q$ implies that in fact $f$ is a strong realization of $L_f$.
\end{proof}

\noindent \textbf{Acknowledgements:} 
The authors are immensely grateful to Professor Osamu Saeki, Director of Institute of Mathematics for Industry at Kyushu University, for discussions and encouragements, and to Benjamin Bode for the comments and recommendations. Part of this paper was developed during the visit of second author at the Kyushu University, Japan, supported by  grants 2019/11415-3 and 2017/25902-8, São Paulo Research Foundation (FAPESP).
\bibliographystyle{amsplain}
  \renewcommand{\refname}{ \large R\normalsize  EFERENCES}
\bibliography{sample}

\providecommand{\bysame}{\leavevmode\hbox to3em{\hrulefill}\thinspace}
\providecommand{\MR}{\relax\ifhmode\unskip\space\fi MR }
\providecommand{\MRhref}[2]{%
  \href{http://www.ams.org/mathscinet-getitem?mr=#1}{#2}
}
\providecommand{\href}[2]{#2}
\begin{thebibliography}{10}

\bibitem{Akbulut_King1981}
S.~Akbulut and H.~King, \emph{All knots are algebraic}, Comment. Math. Helv.
  \textbf{56} (1981), no.~3, 339--351.

\bibitem{Benedetti_Shiota1998}
R.~Benedetti and M.~Shiota, \emph{On real algebraic links in {$S^3$}}, Boll.
  Unione Mat. Ital. Sez. B Artic. Ric. Mat. (8) \textbf{1} (1998), no.~3,
  585--609.

\bibitem{Blanloe_Oka2015}
V.~Blanl{\oe}il and M.~Oka, \emph{Topology of strongly polar weighted
  homogeneous links}, SUT J. Math. \textbf{51} (2015), no.~1, 119--128.

\bibitem{Bode2018}
B.~Bode, \emph{Knotted fields and real algebraic links}, Ph.D. thesis, School
  of Physics, University of Bristol, 2018.

\bibitem{Bode2019}
\bysame, \emph{Constructing links of isolated singularities of polynomials
  {$\Bbb R^4\to\Bbb R^2$}}, J. Knot Theory Ramifications \textbf{28} (2019),
  no.~1, 1950009, 21.

\bibitem{Bode2020}
\bysame, \emph{Real algebraic links in {$S^3$} and braid group actions on the
  set of {$n$}-adic integers}, J. Knot Theory Ramifications \textbf{29} (2020),
  no.~6, 2050039, 44.

\bibitem{Bode2020_arXiv}
\bysame, \emph{Twisting and satellite operations on p-fibered braids}, 2020.

\bibitem{Dennis_Bode2019}
B.~Bode and M.~R. Dennis, \emph{Constructing a polynomial whose nodal set is
  any prescribed knot or link}, J. Knot Theory Ramifications \textbf{28}
  (2019), no.~1, 1850082, 31.

\bibitem{Brieskorn_Knorrer1986}
E.~Brieskorn and H.~Kn{\"o}rrer, \emph{Local investigations}, pp.~325--575,
  Birkh{\"a}user Basel, Basel, 1986.

\bibitem{Dennis_Bode2017}
M.~R. Dennis and B.~Bode, \emph{Constructing a polynomial whose nodal set is
  the three-twist knot {$5_2$}}, J. Phys. A \textbf{50} (2017), no.~26, 265204,
  18.

\bibitem{Garcia_Lenarcik_Ploski2007}
E.~R. Garc\'{\i}a~Barroso, A.~Lenarcik, and A.~P{\l}oski,
  \emph{Characterization of non-degenerate plane curve singularities}, Univ.
  Iagel. Acta Math. \textbf{1298} (2007), no.~45, 27--36.

\bibitem{Inaba_Kawashima_Oka2018}
K.~Inaba, M.~Kawashima, and M.~Oka, \emph{Topology of mixed hypersurfaces of
  cyclic type}, J. Math. Soc. Japan \textbf{70} (2018), no.~1, 387--402.

\bibitem{Kouchnirenko1976}
A.~G. Kouchnirenko, \emph{Poly\`edres de {N}ewton et nombres de {M}ilnor},
  Invent. Math. \textbf{32} (1976), no.~1, 1--31.

\bibitem{Looijenga1971}
E.~Looijenga, \emph{A note on polynomial isolated singularities}, Nederl. Akad.
  Wetensch. Proc. Ser. A {\bf 74}=Indag. Math. \textbf{33} (1971), 418--421.

\bibitem{Milnor1968}
J.~Milnor, \emph{Singular points of complex hypersurfaces}, Annals of
  Mathematics Studies, No. 61, Princeton University Press, Princeton, N.J.;
  University of Tokyo Press, Tokyo, 1968.

\bibitem{Oka2008}
M.~Oka, \emph{Topology of polar weighted homogeneous hypersurfaces}, Kodai
  Math. J. \textbf{31} (2008), no.~2, 163--182.

\bibitem{Oka2010}
\bysame, \emph{Non-degenerate mixed functions}, Kodai Math. J. \textbf{33}
  (2010), no.~1, 1--62.

\bibitem{Oka2015}
\bysame, \emph{Mixed functions of strongly polar weighted homogeneous face
  type}, Singularities in geometry and topology 2011, Adv. Stud. Pure Math.,
  vol.~66, Math. Soc. Japan, Tokyo, 2015, pp.~173--202.

\bibitem{Perron1982}
B.~Perron, \emph{Le n\oe ud ``huit'' est alg\'{e}brique r\'{e}el}, Invent.
  Math. \textbf{65} (1981/82), no.~3, 441--451.

\bibitem{Pichon2005}
A.~Pichon, \emph{Real analytic germs {$f\overline g$} and open-book
  decompositions of the 3-sphere}, Internat. J. Math. \textbf{16} (2005),
  no.~1, 1--12.

\bibitem{Rudolph1987}
L.~Rudolph, \emph{Isolated critical points of mappings from {${\bf R}^4$} to
  {${\bf R}^2$} and a natural splitting of the {M}ilnor number of a classical
  fibered link. {I}. {B}asic theory; examples}, Comment. Math. Helv.
  \textbf{62} (1987), no.~4, 630--645.

\bibitem{Stallings1978}
J.~R. Stallings, \emph{Constructions of fibred knots and links}, Algebraic and
  geometric topology ({P}roc. {S}ympos. {P}ure {M}ath., {S}tanford {U}niv.,
  {S}tanford, {C}alif., 1976), {P}art 2, Proc. Sympos. Pure Math., XXXII, Amer.
  Math. Soc., Providence, R.I., 1978, pp.~55--60.

\bibitem{Zariski1932}
O.~Zariski, \emph{On the {T}opology of {A}lgebroid {S}ingularities}, Amer. J.
  Math. \textbf{54} (1932), no.~3, 453--465.

\end{thebibliography}
\Addresses
\end{document}